\documentclass[11pt]{amsart}
\usepackage{amssymb,amsmath,amsthm}
\usepackage{epsfig}
\newtheorem{theorem}{Theorem}[section]
\newtheorem{lemma}[theorem]{Lemma}

\newtheorem{definition}[theorem]{Definition}
\newtheorem{conjecture}[theorem]{Conjecture}
\newtheorem{proposition}[theorem]{Proposition}
\newtheorem{corollary}[theorem]{Corollary}

\newtheorem{example}[theorem]{Example}
\long\def\symbolfootnote[#1]#2{\begingroup%
\def\thefootnote{\fnsymbol{footnote}}\footnote[#1]{#2}\endgroup}

\newcommand\Z{{\mathbb Z}}

\begin{document}

\title{Homology of dihedral quandles}
\author{Maciej Niebrzydowski and J\'ozef H. Przytycki}
\date{November 25, 2006}
\subjclass{Primary 55N35; Secondary 18G60, 57M25}
\keywords{dihedral quandle, Alexander quandle, rack, quandle homology, homological operations}

\thispagestyle{empty}

\begin{abstract}
We solve the conjecture by R. Fenn, C. Rourke and B. Sanderson 
that the rack homology of dihedral quandles satisfies 
$H_3^R(R_p) = \Z \oplus \Z_p$ for $p$ odd prime 
(\cite{Oht}, Conjecture 5.12). We also show that  
$H_n^R(R_p)$ contains $\Z_p$ for $n\geq 3$. Furthermore, we show that  
the torsion of $H_n^R(R_3)$ is annihilated by $3$.  
We also prove that the quandle homology $H_4^Q(R_p)$ contains $\Z_p$ for $p$ odd prime.
We conjecture that for $n>1$ quandle homology satisfies:
$H_n^Q(R_p) = \Z_p^{f_n}$, where $f_n$ are ``delayed" Fibonacci numbers,
that is, $f_n = f_{n-1} + f_{n-3}$ and $f(1)=f(2)=0, f(3)=1$.
Our paper is the first step in approaching this conjecture.
\end{abstract}

\maketitle

\tableofcontents

\section{Definitions and Preliminary facts}  \label{1}
In 1942, Mituhisa Takasaki introduced the notion of kei \cite{Tak}
as an abstraction of the notion of symmetric transformation. 
He analysed, in particular, finite keis corresponding to symmetries 
of regular polygons (today they are called dihedral quandles).
Kei is an abstract algebra $(Q,*)$
with a binary operation $*: Q\times Q \to Q$ satisfying the conditions:
\begin{enumerate}
\item [(i)] $a*a=a$, for any $a\in Q$;
\item [(ii)] $(a*b)*b=a$, for any $a$, $b\in Q$;
\item [(iii)]$(a*b)*c= (a*c)*(b*c)$, for any $a$, $b$, $c\in Q$ (the right distributivity property).
\end{enumerate}
If the second condition is relaxed to the following condition:\\
(ii') for every $b\in Q$ the map $*_b: Q \to Q$, given by $*_b(a)=a*b$, is
a bijection,
then $(Q,*)$ is called a quandle (the name coined by Joyce \cite{Joy}).\\
A kei is an involutive ($(*_b)^2=Id$) quandle. If the first condition is
omitted, the abstract algebra we obtain is called a rack (J. H. Conway 
in his 1959 correspondence  with G.~C.~Wraith
suggested {\it wrack} which was changed to rack\footnote{From \cite{F-R}:
The earliest work on racks (known to us) is due to Conway and Wraith 
[correspondence, 1959]... They used the name {\it wrack} for the concept 
and we have adopted this name, not merely because it is the oldest name, 
but also because it is a simple English word which (to our knowledge) has 
no other mathematical meaning. We have however chosen the more common 
spelling. Rack is used in the same sense as in the phrase ``rack and ruin".
The context of Conway and Wraith's work is the conjugacy operation in a group 
and they regarded a rack as the wreckage of a group left behind after the 
group operation is discarded and only the notion of conjugacy remains.} 
in \cite{F-R}).
 
Takasaki was considering keis associated to abelian groups, that is,
$kei(G)$ satisfies $a*b= 2b-a$, for $a$, $b\in G$. The $k$-dihedral quandle, denoted $R_k$, is
isomorphic to $kei(\Z_k)$. Generalization of a kei of abelian group 
to Alexander quandle of $Z[t, t^{-1}]$-module 
can be traced back to R.~Fox and his suggestion that arcs of the link diagram can be colored using polynomials.
The Alexander quandle $Alex(M)$ can be 
associated to every $Z[t, t^{-1}]$ module $M$, by taking the quandle operation to be $a*b= (1-t)a +tb$ for any $a$, $b\in M$.

Rack homology and homotopy theory were first defined and 
studied in \cite{FRS}, and a modification to quandle homology theory 
was given in \cite{CJKLS} to define knot invariants in a state-sum form (so-called cocycle knot invariants).

We recall here the definition of rack, degenerate and quandle homology
after \cite{CKS}.
\begin{definition} 
\begin{enumerate}
\item[(i)]
For a given rack $X$, let $C^R_n(X)$ be the free abelian 
group generated by $n$-tuples $(x_1,x_2,\ldots,x_n)$ of elements of $X$; 
in other words, $C^R_n(X) = {\Z}X^n = ({\Z}X)^{\otimes n}$.\\ 
Define a boundary   
homomorphism $\partial: C^R_n(X) \to C^R_{n-1}(X)$ by:
$$\partial(x_1,x_2,\ldots,x_n) = $$ 
$$\sum_{i=2}^n (-1)^i((x_1,\ldots,x_{i-1},x_{i+1},\ldots, 
x_n) - (x_1*x_i,x_2*x_i,\ldots,x_{i-1}*x_i,x_{i+1},\ldots,x_n)).$$ 
$(C^R_*(X),\partial)$ is called a rack chain complex of $X$.
\item[(ii)] Assume that $X$ is a quandle. Then we have a subchain 
complex $C^D_n(X) \subset C^R_n(X)$, generated by $n$-tuples $(x_1,\ldots,x_n)$ 
with $x_{i+1}=x_i$ for some $i$. The subchain complex $(C^D_n(X),\partial)$ 
is called a degenerated chain complex of a quandle $X$.
\item[(iii)] The quotient chain complex $C^Q_n(X)=C^R_n(X)/C^D_n(X)$ is 
called the quandle chain complex. We have the short exact sequence 
of chain complexes:
$$ 0 \to C^D_n(X) \to C^R_n(X) \to C^Q_n(X)\to 0.$$
\item[(iv)] The homology of rack, degenerate, and quandle chain complexes 
are called rack, degenerate, and quandle homology, respectively. We have 
the long exact sequence of homology of quandles:
$$ \ldots\to H^D_n(X) \to H^R_n(X) \to H^Q_n(X)\to H^D_{n-1}(X) \to\ldots$$
\item[(v)] For an abelian group $G$, define the chain complex\\ $C^Q_*(X; G)=C^Q_*\otimes G$, with $\partial=\partial\otimes id$.
The groups of cycles and boundaries are denoted respectively by $ker(\partial)=Z^Q_n(X; G)\subset C^Q_n(X; G)$
and $Im(\partial)=B^Q_n(X; G)\subset C^Q_n(X; G)$.
The $n$th quandle homology group of a quandle $X$ with coefficient group $G$ is defined as
$$H^Q_n(X; G)=H_n(C^Q_*(X; G))=Z^Q_n(X; G)/B^Q_n(X; G).$$
\end{enumerate}
\end{definition}
R.~Litherland and S.~Nelson \cite{L-N} 
proved that the short exact sequence from (iii) splits respecting the
boundary maps.\\ 
The splitting map, $\alpha: C^Q_n(X) \to C^R_n(X)$, is given by:
$$\alpha(x_1,x_2,x_3,\ldots,x_n) = (x_1,x_2-x_1,x_3-x_2,\ldots,x_n-x_{n-1}),$$
where, in our notation, $(x_1,x_2-x_1,x_3-x_2,\ldots,x_n - x_{n-1}) = 
x_1\otimes (x_2-x_1)\otimes (x_3-x_2)\otimes \cdots \otimes 
(x_n-x_{n-1}) \in C^R_n(X)$.
In particular, $\alpha$ is a chain complex monomorphism and 
$H^R_n(X) = H^D_n(X) \oplus \alpha_*(H^Q_n(X))$.
 
Free part of homology of quandles ($free(H_*(X))$) was computed in
\cite{L-N,E-G} (lower bounds for Betti numbers were given in \cite{CJKS}). In particular, it was shown there that for dihedral quandle 
$R_k$, k odd, we have: 
\begin{theorem}\label{Theorem 1.2}
\begin{enumerate}
\item[(i)] \cite{L-N,E-G}
$$free(H^R_n(R_k))=\Z$$
\item[(ii)] \cite{L-N,E-G}
\begin{displaymath}
free(H^Q_n(R_k))=
\left\{
\begin{array}{ll}
\Z & for \ n=1 \\
0 & otherwise
\end{array}
\right.
\end{displaymath}
\end{enumerate}
\end{theorem} 

Useful information concerning torsion of homology of 
racks and quandles was obtained in \cite{L-N,Moc}. 
In particular, it was shown that:
\begin{theorem}\label{Theorem 1.3}
\begin{enumerate}
\item[(i)] \cite{L-N,Moc}
Torsion of $H^R_n(R_k)$ is annihilated by $k^n$.  
\item[(ii)] \cite{Moc} The third cohomology $H^3_Q(R_p;\Z_p)=\Z_p$ for 
$p$ odd prime. 
\item[(iii)] \cite{L-N}\ 
For a quandle $X$, let $\mathcal{O}_X$ denote the set of orbits of $X$ with 
respect to the action of $X$ on itself by the right multiplication (e.g. $|\mathcal{O}_{R_k}|$ is 
$1$ for odd $k$ and $2$ for an even $k$). Then\\
$$ H_2^R(X) \cong H_2^Q(X) \oplus \Z\mathcal{O}_X,$$
$$ H_3^R(X) \cong H_3^Q(X) \oplus H_2^Q(X) \oplus \Z\mathcal{O}_X^2.$$
In particular, $H_3^R(R_k) \cong H_3^Q(R_k) \oplus \Z$ for $k$ odd.
\end{enumerate}
\end{theorem}
We devote this paper to computation of rack and quandle homology 
of odd dihedral quandles.
One of the results of this paper is the solution of the conjecture by Fenn, Rourke and Sanderson, 
listed as Conjecture 5.12 in Ohtsuki's problem list \cite{Oht}.
\begin{conjecture}[R. Fenn, C. Rourke, B. Sanderson]  \ \\  \label{conjRFS}
$H_3^R(R_p) \cong \Z \oplus \Z_p$ for $p$ odd prime.
\end{conjecture}
We prove Conjecture \ref{conjRFS} in Section 3.  We also prove the generalization 
of Theorem 1.3(i) for $k=3$, namely, we show that $3$ annihilates the torsion of 
$H^R_n(R_3)$ for any $n$. We also propose the general conjecture on the structure 
of quandle homology $H^Q_n(R_k)$.
\begin{conjecture} \label{structure}
For a prime number $p$, $tor H^Q_n(R_p) = Z_p^{f_n}$, where 
$\{f_n\}$ are ``delayed" Fibonacci numbers,
that is, $$f_n = f_{n-1} + f_{n-3}, \ \textrm{and}\ f(1)=f(2)=0, f(3)=1.$$
\end{conjecture}
We verified the conjecture by GAP \cite{GAP} calculation for $p=3$ and 
$n\leq 12$, $p=5$ 
and $n\leq 6$, and for $p=7$ and $n\leq 4$.\\ 
We also computed\footnote{ The cases $p=3$, $n=12$, $p=5$, $n=6$ and $p=9$, $n=4$ were computed with 
the help of A.~Shumakovitch and Norbert A'Campo, who provided computer power of his supercomputer located at Basel University, Switzerland.} 
that $H_3^Q(R_9)= H_4^Q(R_9)= \Z_9$, 
so there is a possibility that 
Conjecture \ref{structure} may hold for some non-prime odd $k$ 
(e.g., powers of odd prime numbers). One cannot also exclude that 
it holds for any odd $k$. 

\section{Homological operations}\label{Section 2} 

In this section we construct two chain maps raising index of chain groups 
of a rack (and quandle) by one and two, respectively. We use these maps 
in Sections 3 and 4 to prove our main results on homology of racks and 
quandles.
We expect to find one more homology operation, from $H_n(R_k)$ to
$H_{n+4}(R_k)$, but it is still an open problem.
In this paper, we are mostly interested in homology of dihedral quandles, 
but sometimes we formulate more general results, if it doesn't make proofs much longer.
  
In general, our homological operations involve
the group homomorphism
$h_u: C^R_{n}(X) \to C^R_{n+j}(X)$, for $u\in \Z X^j$,
given by
$h_u (w) = (w,u)$. This map usually is not a chain map, unless
we choose a special $u$ and/or consider quandles satisfying some
special conditions (e.g., as in Lemma 2.5).\\
The first homological operation is related to the group homomorphism\\ 
$h_a\colon C^R_{n}(X) \to C^R_{n+1}(X)$, given by
$h_a(w)=(w,a)$, for any $a\in X$, and $w\in X^n$. This map is not a chain map, so we 
need to symmetrize it with respect to another map $*_a: C^R_{n}(X) \to 
C^R_{n}(X)$ given by $*_a(w)=w*a$, for any $w\in X^n$, or more precisely, 
$*_a(x_1,\ldots,x_n) = (x_1,\ldots,x_n)*a= (x_1*a,\ldots,x_n*a)$ (the map $*_a$ for 
$n=1$ is exactly the map used in Condition (ii') of a quandle). In other words, we consider a function $h'_a = h_a + *_ah_a$. 
The basic properties of these maps are described in the following proposition:
\begin{proposition}\label{Proposition 2.1}\ \
\begin{enumerate}
\item[(i)] For any rack $X$, and $a\in X$, the map $*_a$ is a chain 
map chain homotopic to the identity.
\item[(ii)] If $X$ is a kei (involutive quandle), \label{keichainmap}
then\footnote{To be more
general, we have to assume that $a$ satisfies the ``n-condition",
that is, $x*a*a*\ldots*a=x$ for any $x$, and take $h'_a = h_a + *_ah_a + *_a*_ah_a + \ldots
+ *_a(\ldots(*_ah_a))$.
For example, for quandle $S_4$, $h'_a = h_a + *_ah_a + *_a*_ah_a$. Then
$h'_a$ is a chain map. Notice that if $X$ is finite, then any $a \in X$
satisfies $n$ condition for some $n$.}
$h'_a = h_a + *_ah_a$ is a chain map.
\item[(iii)] If $X$ is a quandle, then $h_a(C^D_{n}(X)) \subset C^D_{n+1}(X)$, 
 $*_a(C^D_{n}(X)) \subset C^D_{n}(X)$.
Therefore, the maps $h_a$, $h'_a: C^Q_{n}(X) \to C^Q_{n+1}(X)$ are well defined.
\item[(iv)] If $a$ and $b$ are in the same orbit of $X$, then 
$h'_a$ and $h'_b$ induce the same map on homology, that is, 
$(h'_a)_*=(h'_b)_*: H^R_{n}(X) \to H^R_{n+1}(X)$.
\end{enumerate}
\end{proposition}
\begin{proof}(i) $*_a$ is a chain map, because\footnote{
We use a standard convention for products in non-associative algebras, called the left normed convention, that is, whenever parentheses are omitted in a product of elements $a_1$, $a_2,\ldots,$ $a_n$ of $Q$ then $a_1*a_2*\ldots *a_n=((\ldots ((a_1*a_2)*a_3)*\ldots)*a_{n-1})*a_n$ (left association), for example, $a*b*c=(a*b)*c$.}
{\small
$$d(*_a(x_1,\ldots,x_n))= d(x_1*a,\ldots,x_n*a)= 
\sum_{i=2}^k (-1)^i((x_1*a,\ldots,x_{i-1}*a,x_{i+1}*a,\ldots,x_n*a)$$ 
$$-(x_1*a*(x_i*a),\ldots,x_{i-1}*a*(x_i*a),x_{i+1}*a,\ldots,x_n*a))=$$
$$\sum_{i=2}^k (-1)^i((x_1*a,\ldots,x_{i-1}*a,x_{i+1}*a,\ldots,x_n*a)-
(x_1*x_i*a,\ldots,x_{i-1}*x_i*a,x_{i+1}*a,\ldots,x_n*a))$$
$$=*_a(\sum_{i=2}^k (-1)^i((x_1,\ldots,x_{i-1},x_{i+1},\ldots,x_n) -
(x_1*x_i,\ldots,x_{i-1}*x_i,x_{i+1},\ldots,x_n)))=$$
$$*_a(d(x_1,\ldots,x_n)).$$}
The homomorphism $(-1)^{n+1}h_a\colon C^R_{n}(X) \to C^R_{n+1}(X)$ is 
a chain homotopy between $Id$ and $*_a$ chain maps. Namely:
$$d((-1)^{n+1}h_a(x_1,\ldots,x_n))=$$ 
$$(-1)^{n+1}((d(x_1,\ldots,x_n),a) +
(-1)^{n+1}((x_1,\ldots,x_n)- (x_1*a,\ldots,x_n*a)))=$$
$$-(-1)^n h_a(d(x_1,\ldots,x_n)) + (Id-*_a)(x_1,\ldots,x_n),$$ as needed.\\
(ii) $dh'_a =d(h_a +*_ah_a)= dh_a + (dh_a)*a =h_ad +(-1)^{n+1}(Id -*_a) + 
(h_ad +(-1)^{n+1}(Id -*_a))*a = h_ad + (h_a*a)d = h'_ad$.\\
(iii) It follows from the definition of rack, degenerate and quandle 
chain complex of $X$.\\
(iv) It suffices to consider the case, when there is $x\in X$, 
such that $b=a*x$.\\ 
Notice, that $$*_x(h'_a(w))= (h'_a(w))*x= 
((w+w*a),a)*x= ((w*x+(w*a)*x),b)=$$ 
$$((w*x+(w*x)*b),b)=h'_b(w*x)=h'_b(*_x(w)).$$
On the other hand, if $w$ is a cycle then, by (i), $*_x(w)$ 
is a homologous cycle. Therefore, $h'_b(*_x(w))$ and $h'_b(w)$ are 
homologous by (ii). 
Similarly, $*_x(h'_a(w))$ and $h'_a(w)$ are homologous. Therefore, $h'_b(w)$ 
and $h'_a(w)$ are homologous.
\end{proof}

Next, we prove that under certain assumptions that hold for odd dihedral quandles,
the map $(h'_a)_*$ is a monomorphism on rack homology. 
\begin{definition}\label{Definition 2.2} 
We say that a rack $X$ satisfies:
\begin{itemize}
\item[--]the property (1), if for any elements $x$, $y$, $a\in X$,
$*_x(a) = *_y(a)$ implies that the maps $*_x$ and $*_y$
are equal (i.e., $z*x = z*y$, for any $z\in X$);
\item[--]the property (2), if for any elements $x$, $y$, $x'$, $y'$, $a\in X$, from\\
$ a*x*x' = a*y*y'$ follows that the maps  $*_{x'}*_x$ and 
$*_{y'}*_y $ are equal (i.e., $z*x'*x = z*y'*y$, for any $z\in X$);
\item[--]the quasigroup property, if for any $a$, $b\in X$, the equation $a*x=b$ has exactly one solution.
\end{itemize}
\end{definition}
An example of racks satisfying properties $(1)$ and $(2)$ is algebraically connected (i.e., with one orbit) Alexander 
quandles.
Indeed, in such case we can write $*_x(a) = *_y(a)$ as $(1-t)x + ta = (1-t)y + ta,$
so $(1-t)x = (1-t)y$, and that implies property $(1)$. 
 
Similarly, $ a*x*x' = a*y*y'$ is equivalent to 
$$(1-t)(x'+tx) + t^2 a = (1-t)(y'+ty) + t^2 a,$$ and further to 
$$(1-t)(x'+tx) = (1-t)(y'+ty),$$ 
from which property $(2)$ follows.

The quasigroup property is stronger than property $(1)$ and an 
Alexander quandle  $M$ possesses this property
only if $(1-t)$ does not annihilate non-zero elements
in the module $M$, and if division by $1-t$ is possible. 
Odd dihedral quandles satisfy all above properties.

\begin{theorem}\label{Theorem 2.3}
Let a kei $X$ be a quasigroup satisfying property (2) of Definition 2.2. 
Then
$(h'_a)_*: H^R_{n}(X) \to H^R_{n+1}(X)$ has the following property.\\
There is a map $\bar h'_*: H^R_{n+1}(X) \to H^R_{n}(X)$ such that
$\bar h'_* (h'_a)_*=4Id$. 
In particular, if $H^R_{n}(X)$ has no element
of order $2$, then $(h'_a)_*$ is a monomorphism.
\end{theorem}

\begin{proof}
From the Proposition \ref{keichainmap}(ii), we know that $h'_a = h_a + *_ah_a$ is a chain map.
Let $\bar h_a : C^R_{n+1}(X) \to C^R_n(X)$ be given by
$$\bar h_a(x_1,\ldots,x_n,x_{n+1}) = (x'_1,\ldots,x'_n),$$ 
where $(x'_1,\ldots,x'_n)$ is uniquely determined by the equality (with unique $x$)
$$(x_1,\ldots,x_n,x_{n+1})*x = (x'_1,\ldots,x'_n,a).$$
We will show that the map $\bar h'_a = \bar h_a + *_a \bar h_a$ is a chain map.
We have
$$d(\bar h'_a(x_1,\ldots,x_n,x_{n+1})) = d((x_1,\ldots,x_n)*x + (x_1,\ldots,x_n)*x*a)=$$
$$(d(x_1,\ldots,x_n))*x  + (d(x_1,\ldots,x_n))*x*a.$$
We also have
$$\bar h'_a(d(x_1,\ldots,x_n,x_{n+1})) = \bar h'_a(d(x_1,\ldots,x_n),x_{n+1}) +$$
$$(-1)^{n+1}(\bar h_a(x_1,\ldots,x_n)- \bar h_a((x_1,\ldots,x_n)*x_{n+1}) +$$
$$(\bar h_a(x_1,\ldots,x_n))*a - (\bar h_a((x_1,\ldots,x_n)*x_{n+1}))*a).$$
To show that  $\bar h'_a$ is a chain map, we need equalities:
$$\bar h_a(x_1,\ldots,x_n) = (\bar h_a((x_1,\ldots,x_n)*x_{n+1}))*a,$$
and
$$(\bar h_a(x_1,\ldots,x_n))*a = \bar h_a((x_1,\ldots,x_n)*x_{n+1}).$$
Notice, that these equalities are equivalent when $X$ is a kei.\\
In order to prove them, we need all assumptions stated in the theorem.
Let $x$ and $y$ be such that $x_n*x=a$ and $x_n*x_{n+1}*y=a$. 
Therefore, $x_n*x*a = x_n*x_{n+1}*y$.
It follows, by property $(2)$, that 
the maps $*_a*_x$ and $*_y*_{x_{n+1}}$ are equal on the 
entire $X$, and the required equalities follow.\\
Once we know that $\bar h'_a$ is a chain map, we can finish the proof by noticing that
$$\bar h'_a h'_a (w) = w + w*a + w*a + w*a*a = 2(w +w*a).$$
Therefore, on homology, $(\bar h'_a h'_a)_* =4Id$.\\
From the Theorem \ref{Theorem 1.3}(i) follows that, for $k$ odd, $H^R_n(R_k)$ can only have odd torsion.
Therefore, for $k$ odd, $(h'_a)_*: H^R_{n}(R_k) \to H^R_{n+1}(R_k)$ is a monomorphism.
\end{proof}

\begin{corollary}\label{Corollary 2.4} 
Let $X$ be a kei satisfying the conditions of Theorem \ref{Theorem 2.3}.\\
(i) If $tor H_n^R(X)$ is annihilated by $N$, then 
$tor H_{n+1}^R(X)$ is annihilated 
by $4|X|N$. \\
(ii) If $tor H_n^R(R_k)$, $k$ odd, is annihilated by $N$ 
then $tor H_{n+1}^R(R_k)$ is annihilated by $kN$.\\
(iii) For $k$ odd, $tor H_3^R(R_k)$ is annihilated by $k$.
\end{corollary}
\begin{proof}(i) Consider the (symmetrizer) chain map 
$s_{n+1}: C_{n+1}(X) \to C_{n+1}(X)$ given by
$$s_{n+1}(u) = \sum_{x\in X}(u*x).$$ By Proposition \ref{Proposition 2.1}(i) the map 
$(s_{n+1})_*: H^R_{n+1}(X) \to H^R_{n+1}(X)$ is equal to  $|X|Id$. 
Next, consider the composition $$s_{n+1}h'_a\bar h'_a : C_{n+1}(X) \to C_{n+1}(X).$$
To find $s_{n+1}h'_a\bar h'_a(w,x_{n+1})$, we assume that $x_{n+1}*x=a$ (so that we also know that $x_{n+1}=a*x$).
Then: 
$$s_{n+1}h'_a\bar h'_a(w,x_{n+1})= s_{n+1}h'_a(w*x + w*x*a)= 
s_{n+1}(2((w*x,a) + (w*x*a,a)))=$$ 
$$2\sum_{y\in X}((w*x*y,a*y) + (w*x*a*y,a*y))=
2\sum_{y\in X}((w,x_{n+1})*y*a + (w,x_{n+1})*y).$$
In order to prove the last equality, we need to show that:
\begin{itemize}
\item[(1)] $\sum_{y\in X}(w*x*y,a*y)=\sum_{y\in X}(w,x_{n+1})*y*a$;
\item[(2)] $\sum_{y\in X}(w*x*a*y,a*y)=\sum_{y\in X}(w,x_{n+1})*y$.
\end{itemize}
The left hand side of (1) is equal to
$$\sum_{y\in X}(w*x*y,a*y)=\sum_{y\in X}(w*x,a)*y = \sum_{y\in X}(w,x_{n+1})*x*y.$$  
Now, the set
$\{x_{n+1}*x*y\colon\ y\in X\}$ is the set of all elements of $X$,
as is the set 
$\{x_{n+1}*y*a\colon\ y\in X\}$.
Furthermore, by Property (2) of Definition \ref{Definition 2.2},
if
$x_{n+1}*x*y= x_{n+1}*y'*a$ for some $y$ and $y'$, then
$w*x*y= w*y'*a$ for any $w\in C_n(X)$.
Thus, the equality (1) holds.\\
To prove $(2)$, we notice that:
$$\sum_{y\in X}(w*x*a*y,a*y)=\sum_{y\in X}(w*x,a)*(y*a)*a=$$
$$\sum_{y\in X}(w,x_{n+1})*x*(y*a)*a=\sum_{y\in X}(w,x_{n+1})*x*y*a.$$
By a similar argument as in $(1)$, for every $y\in X$ there exists $y'\in X$ such that
$x_{n+1}*x*y*a=x_{n+1}*y'$, i.e., $x_{n+1}*x*y=x_{n+1}*y'*a$.
It follows that $w*x*y=w*y'*a$ (and $w*x*y*a=w*y'$) for any $w\in C_n(X)$, and that finishes the proof of equality (2).\\
By Proposition \ref{Proposition 2.1}(i) 
the map $(s_{n+1}h'_a\bar h')_* = 4|X|Id$. Therefore, if $tor H_n^R(X)$ is 
annihilated by $N$, then $tor H_{n+1}^R(X)$ is annihilated
by $4|X|N$. \\
(ii) By Theorem \ref{Theorem 1.3}(i), $H_n^R(R_k)$, $k$ odd, has no even torsion. Thus, 
if $tor H_n^R(X)$ is annihilated by $N$, then $tor H_{n+1}^R(X)$ is annihilated
by $kN$.\\ 
(iii) As given in [CJKS] (Example 3.10), $tor H_2^R(R_k) = 0$ for $k$ odd. Thus, $tor H_3^R(R_k)$ is annihilated by $k$.
\end{proof}

The second homological operation, of degree two, is defined on somehow 
restricted class of racks (and quandles), but including dihedral quandles 
and Burnside keis (see footnote 5).
Let $X$ be any rack with 
elements $x_0$, $x_1,\ldots, x_{k-1}$, satisfying the ``cyclic" Fibonacci 
relation\footnote{If $X$ is a kei, then by eliminating all generators except $x_0$ and $x_1$, we 
obtain relations 
$x_0=((((\ldots x_0)*x_1)*x_0)*x_1)$ and $x_1=((((\ldots x_1)*x_0)*x_1)*x_0)$, with
$k$ letters on the right hand side of each equation. They were studied in \cite{N-P} 
and led to {\it Burnside Keis} which are important in analysis 
of rational moves on links and tangles.} modulo $k$.
That is, $x_2=x_0*x_1$, $x_3=x_1*x_2, \ldots, x_0=x_{k-2}*x_{k-1}$
and $x_1=x_{k-1}*x_{0}$.\\ 
Let $s=s(x_0,x_1)= \sum_{i=0}^{k-1}(x_i,x_{i+1})$, 
where indices are taken modulo $k$.\\
Then we can use $s$ to define a chain map as follows.
\begin{lemma}\label{Lemma 2.6}\
\begin{enumerate}
\item[(i)]
The homomorphism $h_s:C^R_n(X) \to C^R_{n+2}(X)$ given by 
$h_s(w) = (w,s)$, is a chain map. Recall that $C^R_{n+2}(X)= \Z X^{n+2}= 
\Z X^{n}\otimes \Z X^2$ so $(w,s)$ is another notation for $w\otimes s$.
\item[(ii)] If $X=R_k$, then we can take $x_i=i$, that 
is, $$s=s(0,1)= \sum_{i=0}^{k-1}(i,i+1).$$ 
For an odd $k$, we can choose any $0<j<k$, and take
$x_0=0$, $x_1=j$, and generally, $x_i=ij$. Let $u_j\in C_2^R(R_k)$ be defined 
by $u_j=s(0,j)$. 
Then $h_{u_j}: C^R_n(R_k) \to C^R_{n+2}(R_k)$ given by $h_{u_j}(w) = (w,u_j)$, is a chain map.
\end{enumerate}
\end{lemma}
\begin{proof}
(i) We notice that
{\small
$$d(h_s(w))=d(w,s)=$$
$$(dw,s) + \sum_{i=0}^k(-1)^{n+1}((w,x_{i+1})-(w,x_i)) 
- \sum_{i=0}^k(-1)^{n+1}((w*x_i,x_{i+1})-(w*x_{i+1},x_i*x_{i+1}))=$$
$$h_s(d(w)) -\sum_{i=0}^k(-1)^{n+1}((w*x_i,x_{i+1})-(w*x_{i+1},x_{i+2})) = h_s(d(w)),$$}
as required.\\
(ii) We check 
that $ij*(ij+j)= 2ij+2j-ij = ij+2j$, so the ``cyclic" Fibonacci relation is satisfied and we can use part (i).
\end{proof}

We prove, in the next section, that for $k$ odd prime, $h_s: H^Q_1(R_k) \to H^Q_3(R_k)$ is 
an epimorphism $\Z \to \Z_k$. We conjecture the following.
\begin{conjecture} For an odd prime number $k$, and $n>1$, the 
homomorphism $h_s: H^Q_n(R_k) \to H^Q_{n+2}(R_k)$ is a monomorphism.
\end{conjecture}
\begin{example}

We checked (using GAP) that $ (h_s)_*: H^Q_3(R_3) \to H^Q_5(R_3)$ is an isomorphism, sending
$(0,1,0) + (0,2,1)$ to $(0,1,0,1,2) + (0,1,0,2,0) + (0,2,1,0,1) + (0,2,1,2,0)$.
We also checked that $(h_s)_*: H_n(R_3) \to H_{n+2}(R_3)$ is 
a monomorphism for $n=4,5,6,7$. Finally, we checked that $ (h_s)_*: H_3(R_5) \to H_5(R_5)$ is an 
isomorphism.
\end{example}

\section{The proof of the Fenn-Rourke-Sanderson conjecture}
Recall that $u_j \in C_2^R(R_k)$ is given by 
$u_j = \sum_{i=0}^{k-1} (i,i+j)$ and that\\
$(0,u_j) = \sum_{i=0}^{k-1} (0,i,i+j).$\\ 
The generator of the third quandle cohomology $H^3_Q(R_p;\Z_p)$, given by T.~Mochizuki
in \cite{Moc}, is of the form
$$\theta(x_1,x_2,x_3)= \frac{(x_2-x_1)(2x_3^p-x_2^p -(2x_3-x_2)^p)}{p}.$$
\begin{proposition}\label{Evaluation}
$\theta((0,u_j)) = j^2\theta((0,u_1))= -j^2 $ in $\Z_p$. In particular,
$(0,u_j)$ is a nontrivial cycle in $H^Q_3(R_p)$.
\end{proposition}

\begin{proof}
In the straightforward evaluation of the Mochizuki's $3$-cocycle on $(0,u_j)$, we use the easy to check 
fact that $(p+a)^p \equiv a^p \mod p^2$, for any $a\in \Z_p$, and that
$\theta((0,u_j))= \sum_{i=0}^{p-1} -(i)\frac{(i+2j)^p +i^p - 2(i+j)^p}{p}$.
\end{proof}

To show that $H_3^Q(R_p)=\Z_p$, we use the universal coefficient theorem 
for cohomology (see \cite{Ha}, Theorem 3.2).
\begin{theorem}\label{Univ Coef Cohomology}
If a chain complex $C$ of free abelian groups has homology groups $H_n(C)$, 
then the cohomology groups $H^n(C;G)$ of the cochain complex $Hom(C_n,G)$ 
are determined by split exact sequences
$$ 0 \to Ext(H_{n-1}(C),G) \to H^n(C;G) \to Hom(H_n(C),G) \to 0 .$$
\end{theorem}
We also recall that $Ext(H,G)=0$ if $H$ is free.

Now we can complete the proof of the Fenn-Rourke-Sanderson conjecture.
First, recall that $H_3^R(R_p) = H^Q_3(R_p) \oplus \Z$, 
as summarized in Theorem 1.3(iii). 
Thus, we should show that $H^Q_3(R_p) = \Z_p$. 
By Mochizuki's result (Theorem 1.3(ii)), $H^3_Q(R_p;\Z_p) = \Z_p$. Thus, 
from the universal coefficient theorem for cohomology we have:
$$\Z_p = H^3_Q(R_p;\Z_p)= Hom(H_3^Q(R_p),\Z_p) \oplus Ext(H_2^Q(R_p),\Z_p) = 
Hom(H_3^Q(R_p),\Z_p),$$ as $H_2^Q(R_p)=0.$ 
Therefore, $H_3^Q(R_p)$ is equal 
to $Z$ or $\Z_{p^i}$, for $i>0$. 
But by the Corollary \ref{Corollary 2.4} and Theorem \ref{Theorem 1.2}(ii), $H_3^Q(R_p)$ is annihilated by $p$. 
Therefore, $H_3^Q(R_p)=\Z_p$, what ends the proof.

It follows from the definition of rack and quandle homology, that $H_1^R(R_p)$ and $H_1^Q(R_p)$ are generated 
by $(0) \in C_1(R_p)$. It follows from above proofs that $(0,u_1) = \sum_{i=0}^{p-1} (0,i,i+1)$ is a generator of
$H_3^Q(R_p)$, for p odd prime.
Therefore, we conclude that the map $h_s: H^Q_1(R_p) \to H^Q_3(R_p)$, that we discussed in the previous section, is 
an epimorphism $\Z \to \Z_p$.

\section{$3$ annihilates the torsion of $H^R_n(R_3)$}

In this section we show that $3$ annihilates $tor H^R_n(R_3)$, for $n>1$.
In the proof we use the fact that $R_3$ is commutative, i.e., $$a*b=2b-a=-b-a=2a-b=b*a,$$ for any $a, b\in R_3$.
However, there is a possibility that our proof can be generalized (provided that the Conjecture \ref{con2} is true), therefore, 
we write the proof below in the more general form and we conjecture the following:
\begin{conjecture}
If $q$ is odd prime, and $n>1$, then the torsion 
subgroup of $H^R_n(R_q)$ is annihilated by $q$.
\end{conjecture}

\begin{theorem}
If $q=3$, and $n>1$, then $tor H^R_n(R_q)$ is annihilated by $q$.
\end{theorem}

\begin{proof} 
Let $x=(x_1,\ldots,x_n) \in (R_q)^n$.
We define the chain map\\ $f_n^j\colon C^R_n(R_q)\to C^R_n(R_q)$ by

\begin{displaymath}
f_n^j(x)= \left\{ \begin{array}{ll}
x & \textrm{for } j=0\\
\sum_{y\in R_q} (y,\ldots,y,x_{j+1},x_{j+2},\ldots,x_n) & \textrm{for } 1\leq j\leq n\\
f^n_n(x) & \textrm{for } j>n
\end{array} \right.
\end{displaymath}
It was shown in \cite{L-N} that $f^1_n$ is chain homotopic to $q$ times the identity.
We prove that $f^j_n$ is chain homotopic to $f^{j-1}_n$, for $j\geq 2$, and the chain homotopy is given by the formula
\begin{displaymath}
D_n^j(x)= \left\{ \begin{array}{ll}
\sum_{y\in R_q} (y,\ldots,y,x_j,y,x_{j+1},\ldots,x_n) & \textrm{for } 1\leq j\leq n\\
0 & \textrm{for } j>n
\end{array} \right.
\end{displaymath}
We will use standard notation:
$$\partial^0_{n}(x)=\sum_{i=1}^n(-1)^i(x_1,\ldots,x_{i-1},x_{i+1}\ldots,x_n);$$
$$\partial^1_{n}(x)=\sum_{i=1}^n(-1)^{i+1}(x_1*x_i,\ldots,x_{i-1}*x_i,x_{i+1},\ldots,x_n).$$
We need to confirm the equality
$$\partial_{n+1} D^j_n(x)+D^j_{n-1}\partial_n(x)=\pm(f^j_n(x)-f^{j-1}_n(x)).$$
Let us write a detailed list of summands appearing on the left hand side of this equation:

$$\partial^0_{n+1} D_n^j(x)=$$
\begin{equation}
\sum_{i=1}^{j-1}\sum_{y\in R_q} (-1)^i(y,\ldots,y,x_j,y,x_{j+1},\ldots,x_n) + \label{1a}
\end{equation}
\begin{equation}
(-1)^j\sum_{y\in R_q}(y,\ldots,y,x_{j+1},\ldots,x_n) + \label{2a1}
\end{equation}
\begin{equation}
(-1)^{j+1}\sum_{y\in R_q}(y,\ldots,y,x_j,x_{j+1},\ldots,x_n) + \label{2a2}
\end{equation}
\begin{equation}
(-1)^{j+2}\sum_{y\in R_q}(y,\ldots,y,x_j,y,x_{j+2},\ldots,x_n) + \label{3a}
\end{equation}
\begin{equation}
\sum_{i=j+2}^{n}\sum_{y\in R_q}(-1)^{i+1}(y,\ldots,y,x_j,y,x_{j+1},\ldots,\widehat{x_i},\ldots,x_n)  \label{4a}
\end{equation}
$$\partial^1_{n+1} D_n^j(x)=$$
\begin{equation}
\sum_{i=1}^{j-1}\sum_{y\in R_q}(-1)^{i+1}(y,\ldots,y,x_j,y,x_{j+1},\ldots,x_n)  + \label{1b}
\end{equation}
\begin{equation}
(-1)^{j+1}\sum_{y\in R_q}(y*x_j,\ldots,y*x_j,y,x_{j+1},\ldots,x_n) + \label{2b}
\end{equation}
\begin{equation}
(-1)^{j+2}\sum_{y\in R_q}(y,\ldots,y,x_j*y,x_{j+1},\ldots,x_n) + \label{3b}
\end{equation}
\begin{equation}
(-1)^{j+3}\sum_{y\in R_q}(y*x_{j+1},\ldots,y*x_{j+1},x_j*x_{j+1},y*x_{j+1},x_{j+2},\ldots,x_n) + \label{4b}
\end{equation}
\begin{equation}
\sum_{i=j+2}^{n}\sum_{y\in R_q}(-1)^{i+2}(y*x_i,\ldots,y*x_i,x_j*x_i,y*x_i,x_{j+1}*x_i,\ldots,
x_{i-1}*x_i,x_{i+1},\ldots,x_n) \label{5b}
\end{equation}
$$D_{n-1}^j\partial^0_n=$$
\begin{equation}
\sum_{i=1}^j\sum_{y\in R_q}(-1)^i(y,\ldots,y,x_{j+1},y,x_{j+2},\ldots,x_n) +  \label{1c}
\end{equation}
\begin{equation}
(-1)^{j+1}\sum_{y\in R_q}(y,\ldots,y,x_j,y,x_{j+2},\ldots,x_n) +  \label{2c}
\end{equation}
\begin{equation}
\sum_{i=j+2}^{n}\sum_{y\in R_q}(-1)^i(y,\ldots,y,x_j,y,x_{j+1},\ldots,\widehat{x_i},\ldots,x_n)  \label{3c}
\end{equation}
$$D_{n-1}^j\partial^1_n=$$
\begin{equation}
\sum_{i=1}^j\sum_{y\in R_q}(-1)^{i+1}(y,\ldots,y,x_{j+1},y,x_{j+2},\ldots,x_n) +  \label{1d}
\end{equation}
\begin{equation}
(-1)^{j+2}\sum_{y\in R_q}(y,\ldots,y,x_j*x_{j+1},y,x_{j+2},\ldots,x_n) +  \label{2d}
\end{equation}
\begin{equation}
\sum_{i=j+2}^{n}\sum_{y\in R_q}(-1)^{i+3}(y,\ldots,y,x_j*x_i,y,x_{j+1}*x_i,\ldots,\widehat{x_i},\ldots,x_n)  \label{3d}
\end{equation}
In the above sums, $\widehat{x_i}$ denotes the omission of the element $x_i$.\\
Notice, that we have the following cancellations:
$$(\ref{1a}) + (\ref{1b})=0,$$
$$(\ref{4a}) + (\ref{3c})=0,$$
$$(\ref{3a}) + (\ref{2c})=0,$$
$$(\ref{5b}) + (\ref{3d})=0,$$
$$(\ref{4b}) + (\ref{2d})=0,$$
$$(\ref{1c}) + (\ref{1d})=0.$$
Moreover,
$$(\ref{2a1})+(\ref{2a2})=(-1)^j(f^j_n(x)-f^{j-1}_n(x)).$$
To finish the proof of the chain homotopy between  $f^j_n(x)$ and $f^{j-1}_n(x)$, we need to look closer at
$$(\ref{2b})=(-1)^{j+1}\sum_{y\in R_q}(y*x_j,\ldots,y*x_j,y,x_{j+1},\ldots,x_n),$$
and
$$(\ref{3b})=(-1)^{j+2}\sum_{y\in R_q}(y,\ldots,y,x_j*y,x_{j+1},\ldots,x_n).$$
Notice, that
$$\sum_{y\in R_q}(y*x_j,\ldots,y*x_j,y,x_{j+1},\ldots,x_n)=
\sum_{y\in R_q}(y,\ldots,y,y*x_j,x_{j+1},\ldots,x_n).$$
Therefore, in the case of $R_3$, the commutativity property gives the final cancellation:
$$(\ref{2b}) + (\ref{3b})=0.$$
In the case of prime $q>3$, we could consider maps $g^j_1$, $g^j_2\colon C^R_n(R_q)\to C^R_n(R_q)$ given by:
$$g^j_1(x)=\sum_{y\in R_q}(y,\ldots,y,y*x_j,x_{j+1},\ldots,x_n),$$
and
$$g^j_2(x)=\sum_{y\in R_q}(y,\ldots,y,x_j*y,x_{j+1},\ldots,x_n).$$
These are chain maps, if $q$ is odd and not divisible by $3$. The straightforward proof of this statement uses the fact that, 
for such $q$, the map $r_a\colon R_q\to R_q$ given by $r_a(b)=b*a*b$ is a bijection.\\
We conjecture the following:
\begin{conjecture}\label{con2}
The maps $g^j_1$ and $g^j_2\colon C^R_n(R_q)\to C^R_n(R_q)$ are chain homotopic, for prime $q>3$.
\end{conjecture}
If this conjecture is true, then the maps $f^j_n(x)$ and $f^{j-1}_n(x)$ are chain homotopic for any odd prime $q$ (and $j\geq 2$).
This would allow us to prove the Conjecture 4.1.
For now, we can finish the proof for $q=3$.\\
Recall that $\mathcal{O}(Q)$ denotes the set of orbits of a quandle $Q$. We can consider it as a trivial quandle ($a*b=a$, for any
$a$, $b\in \mathcal{O}(Q)$).\\
Let us define two maps, $\pi\colon C_n(R_q)\to C_n(\mathcal{O}(R_q))$, and $\Psi\colon C_n(\mathcal{O}(R_q))\to C_n(R_q)$ in the following way:
$$\pi(x_1,\ldots,x_n)=(1,\ldots,1),$$
where $1$ denotes the single orbit of the quandle $R_q$;
$$\Psi(1,\ldots,1)=\sum_{y\in R_q} (y,\ldots,y).$$
Notice that both $\Psi$ and $\pi$ are chain maps.
We have $$\Psi\pi(x)=\sum_{y\in R_q} (y,\ldots,y)=f^n_n(x).$$
Therefore, on homology, $\Psi_*\pi_*=q\cdot Id$. It follows that the torsion of $H^R_n(R_q)$ is annihilated by $q$, for $q=3$. 
\end{proof}

\section{$H^Q_4(R_p)$ is nontrivial for $p$ odd prime}
\begin{theorem}
$H^Q_4(R_p)$ contains $Z_p$, for $p$ odd prime.
\end{theorem}
\begin{proof}
As we already know from Section 3,
$$\sum_{i=0}^{p-1} (0,i,i+1)$$ is a generator of $H^Q_3(R_p)$. Furthermore,
$$\alpha(\sum_{i=0}^{p-1} (0,i,i+1))$$ is a generator of $H^R_3(R_p)$, because 
$\alpha_*\colon H^Q_*(R_p)\to H^R_*(R_p)$ is a monomorphism.
$$\Psi=(\alpha(\sum_{i=0}^{p-1} (0,i,i+1)),0) + (\alpha(\sum_{i=0}^{p-1} (0,i,i+1)),0)*0$$ is nontrivial in $H^R_4(R_p)$
(the map $h'_0=h_0+(h_0)*a$ induces monomorphism on homology).\\
We will show that
$$\sum_{i=0}^{p-1} (0,i,i+1,0) + (\sum_{i=0}^{p-1} (0,i,i+1,0))*0$$ is nontrivial in $H^Q_4(R_p)$ by proving that
$$\Phi=\alpha(\sum_{i=0}^{p-1} (0,i,i+1,0) + (\sum_{i=0}^{p-1} (0,i,i+1,0))*0)=$$
$$\alpha(\sum_{i=0}^{p-1} (0,i,i+1,0)) + (\alpha(\sum_{i=0}^{p-1} (0,i,i+1,0)))*0$$ is a cycle homologous to $\Psi$.
That is, we prove that $\Phi-\Psi$ is zero in homology.\\
Notice that 
$$\alpha(0,i,i+1)=(0,i,i+1)-(0,0,i+1)-(0,i,i)+(0,0,i),$$
and therefore,
$$\alpha(\sum_{i=0}^{p-1} (0,i,i+1))=\sum_{i=0}^{p-1}(0,i,i+1)-\sum_{i=0}^{p-1}(0,i,i).$$
Also,
$$\alpha(0,i,i+1,0)=(\alpha(0,i,i+1),0)-(\alpha(0,i,i+1),i+1);$$
thus,
$$\alpha(\sum_{i=0}^{p-1} (0,i,i+1,0))=\sum_{i=0}^{p-1}(0,i,i+1,0)-\sum_{i=0}^{p-1}(0,i,i,0)$$
$$-\sum_{i=0}^{p-1}(0,i,i+1,i+1)+\sum_{i=0}^{p-1}(0,i,i,i+1).$$
It follows that
$$\Phi-\Psi=\sum_{i=0}^{p-1}(0,i,i,i+1) - \sum_{i=0}^{p-1}(0,i,i+1,i+1)
+ (\sum_{i=0}^{p-1}(0,i,i,i+1) - \sum_{i=0}^{p-1}(0,i,i+1,i+1))*0.$$
Notice that $$c=\sum_{i=0}^{p-1}(0,i,i,i+1) - \sum_{i=0}^{p-1}(0,i,i+1,i+1)$$
is a cycle homologous to $$(\sum_{i=0}^{p-1}(0,i,i,i+1) - \sum_{i=0}^{p-1}(0,i,i+1,i+1))*0.$$
Therefore, to show that $(\Phi-\Psi)\sim 0$, it is enough to show that $c\sim 0$
(here, $\sim$ denotes homologous cycles).\\
It follows from the Proposition \ref{Proposition 2.1}(iv) that $$h'_{i+1}(h'_i(h'_i(0))) - h'_{i+1}(h'_{i+1}(h'_i(0))) \sim 0.$$
We have 
{\small $$h'_{i+1}(h'_i(h'_i(0)))=h'_{i+1}(h'_i((0,i)+(0*i,i)))=h'_{i+1}(2((0,i,i)+(0*i,i,i)))=$$
$$2((0,i,i,i+1)+(0*i,i,i,i+1)+(0*(i+1),i+2,i+2,i+1)+(0*i*(i+1),i+2,i+2,i+1)),$$}
and
{\small
$$h'_{i+1}(h'_{i+1}(h'_i(0)))=h'_{i+1}(h'_{i+1}((0,i)+(0*i,i)))=$$
$$h'_{i+1}((0,i,i+1)+(0*i,i,i+1)+(0*(i+1),i+2,i+1)+(0*i*(i+1),i+2,i+1))=$$
$$2((0,i,i+1,i+1)+(0*i,i,i+1,i+1)+(0*(i+1),i+2,i+1,i+1)+(0*i*(i+1),i+2,i+1,i+1)).$$}
From above, and the fact that there can be only odd torsion in homology of $R_k$ with $k$ odd, follows
$$(\sum_{i=0}^{p-1}(0,i,i,i+1)+\sum_{i=0}^{p-1}(2i,i,i,i+1)+\sum_{i=0}^{p-1}(2i,i+1,i+1,i)+
\sum_{i=0}^{p-1}(2,i+1,i+1,i)$$
$$-\sum_{i=0}^{p-1}(0,i,i+1,i+1)-\sum_{i=0}^{p-1}(2i,i,i+1,i+1)-\sum_{i=0}^{p-1}(2i,i+1,i,i)
-\sum_{i=0}^{p-1}(2,i+1,i,i))\sim 0.$$
The above expression decomposes into the following cycles:
\begin{equation}
c=\sum_{i=0}^{p-1}(0,i,i,i+1)-\sum_{i=0}^{p-1}(0,i,i+1,i+1), \label{c1}
\end{equation}
\begin{equation}
\sum_{i=0}^{p-1}(2,i+1,i+1,i)-\sum_{i=0}^{p-1}(2,i+1,i,i),   \label{c2}
\end{equation}
\begin{equation}
\sum_{i=0}^{p-1}(2i,i,i,i+1)-\sum_{i=0}^{p-1}(2i,i+1,i,i),   \label{c3}
\end{equation}
\begin{equation}
\sum_{i=0}^{p-1}(2i,i+1,i+1,i)-\sum_{i=0}^{p-1}(2i,i,i+1,i+1).   \label{c4}
\end{equation}
We will show that all these cycles are homologous to $c$.\\
Let $x$ be such that $2*x=0$. Then
$$(\ref{c2})*x=(\sum_{i=0}^{p-1}(2,i+1,i+1,i)-\sum_{i=0}^{p-1}(2,i+1,i,i))*x=$$
$$\sum_{i=0}^{p-1}(0,2x-(i+1),2x-(i+1),2x-i)-\sum_{i=0}^{p-1}(0,2x-(i+1),2x-i,2x-i)=$$
$$\sum_{i=0}^{p-1}(0,i,i,i+1)-\sum_{i=0}^{p-1}(0,i,i+1,i+1)=(\ref{c1}).$$
Thus, $(\ref{c1})\sim (\ref{c2})$.\\
Now, $(\ref{c1})-(\ref{c3})$ is the following boundary:
{\small
$$\partial(\sum_{i=0}^{p-1}(0,i,i,i,i+1)+\sum_{i=0}^{p-1}(0,i,i+1,i+1,i+1))=$$
$$\sum_{i=0}^{p-1}(0,i,i,i+1)-\sum_{i=0}^{p-1}(2i,i,i,i+1)-\sum_{i=0}^{p-1}(0,i,i,i)
+\sum_{i=0}^{p-1}(2i+2,i+2,i+2,i+2)$$
$$+\sum_{i=0}^{p-1}(0,i+1,i+1,i+1)-\sum_{i=0}^{p-1}(2i,i+1,i+1,i+1)$$
$$-\sum_{i=0}^{p-1}(0,i,i+1,i+1)+\sum_{i=0}^{p-1}(2i+2,i+2,i+1,i+1)=$$
$$\sum_{i=0}^{p-1}(0,i,i,i+1)-\sum_{i=0}^{p-1}(2i,i,i,i+1)
-\sum_{i=0}^{p-1}(0,i,i+1,i+1)+\sum_{i=0}^{p-1}(2i,i+1,i,i)=(\ref{c1})-(\ref{c3}).$$}
Finally, we prove that $(\ref{c2})-(\ref{c4})$ is a boundary.
{\small $$\partial(\sum_{i=0}^{p-1}(2,i+1,i+1,i+1,i)+\sum_{i=0}^{p-1}(2,i+1,i,i,i))=$$
$$\sum_{i=0}^{p-1}(2,i+1,i+1,i)-\sum_{i=0}^{p-1}(2i,i+1,i+1,i)-\sum_{i=0}^{p-1}(2,i+1,i+1,i+1)
+\sum_{i=0}^{p-1}(2i-2,i-1,i-1,i-1)$$
$$+\sum_{i=0}^{p-1}(2,i,i,i)-\sum_{i=0}^{p-1}(2i,i,i,i)-\sum_{i=0}^{p-1}(2,i+1,i,i)
+\sum_{i=0}^{p-1}(2i-2,i-1,i,i)=$$
$$\sum_{i=0}^{p-1}(2,i+1,i+1,i)-\sum_{i=0}^{p-1}(2i,i+1,i+1,i)-\sum_{i=0}^{p-1}(2,i+1,i,i)
+\sum_{i=0}^{p-1}(2i,i,i+1,i+1)=(\ref{c2})-(\ref{c4}).$$}
It follows that
$$4c\sim 0,$$ and this forces $c\sim 0$ (there can be no even torsion).
Thus, $(\Phi-\Psi)\sim 2c \sim 0$, what ends the proof. 
\end{proof}

\section{Future directions}
We have constructed, in Section 2, two homological operations\\
$h'_0: C^R_n(R_k) \to C^R_{n+1}(R_k)$ and
$h_s : C^R_n(R_k) \to C^R_{n+2}(R_k)$.
We have proven that $h'_0$ induces a monomorphism on rack
homology $H^R_n(R_k) \to H^R_{n+1}(R_k)$ for $k$ odd.
We have conjectured that $(h_s)_*: H^R_n(R_k) \to H^R_{n+2}(R_k)$
is also a monomorphism. More generally, we propose:
\begin{conjecture}
$(h'_0)_* \oplus (h_s)_*: H^R_{n+1}(R_k) \oplus H^R_n(R_k) \to H^R_{n+2}(R_k)$
is a monomorphism for $k$ odd and $n>1$.
\end{conjecture}

In the case of quandle homology, the composition $h'_ah'_b$ induces the zero map on
homology for any $a$, $b\in R_k$. Therefore, instead of this map, we consider the composition map
$h'_{(s,0)}= h'_0h_s: C^Q_n(R_k) \to C^Q_{n+3}(R_k)$ and propose:

\begin{conjecture}\label{sprytna}
$(h_s)_* \oplus (h'_{(s,0)})_*: H^Q_{n+1}(R_k) \oplus H^Q_n(R_k) \to H^Q_{n+3}(R_k)$
is a monomorphism for $k$ odd and $n>1$.
\end{conjecture}

We have checked that the conjecture holds for $k=3,5$ and $n=2, 3$ (in these cases,
these maps are isomorphisms) and for $k=3$ and $n=4,5$.

Conjecture \ref{sprytna} is not sufficient to prove our main conjecture (Conjecture 1.5).
To achieve this, we need additional homological operation of degree $4$:
$h_{new}: C^Q_n(R_k) \to C^Q_{n+4}(R_k)$ such that
$$(h_s)_* \oplus (h'_{(s,0)})_* \oplus (h_{new})_*: H^Q_{n+2}(R_k) \oplus H^Q_{n+1}(R_k) \oplus H^Q_n(R_k)
\to H^Q_{n+4}(R_k)$$
is an isomorphism for $k$ odd and $n>1$. It would provide inductive step in a proof of Conjecture 1.5,
as $f_{n+4}=f_{n+3}+f_{n+1} = f_{n+2} + f_{n+1} + f_{n}$. 

\section{Appendix 1}
In this section we demonstrate a generalization of Corollary 2.4.
\begin{theorem}
Let a quandle $Q$ satisfy the quasigroup property (Definition 2.2).
Then, if $N$ annihilates the torsion of $H^R_n(Q)$, then $|Q|N$ annihilates the torsion of
$H_{n+1}^R(Q)$.
\end{theorem}

\begin{proof}
Consider the composition of two homomorphisms $$gf\colon C^R_n(Q) \to  C^R_{n-1}(Q) \to C^R_n(Q), $$ where
$$f(x_1,x_2,\ldots,x_n) = (-1)^n(x_2,\ldots,x_n)$$ and
$$g(x_2,\ldots,x_n) = (-1)^n\sum_{y\in Q}(y,x_2,\ldots,x_n).$$
We check that $f$ and $g$ are chain maps.
{\small
$$fd(x_1,x_2,\ldots,x_n)=$$ $$(-1)^{n-1}\sum_{i=2}^n (-1)^i((x_2,\ldots,x_{i-1},x_{i+1},\ldots,x_n)
-(x_2*x_i,\ldots,x_{i-1}*x_i,x_{i+1},\ldots,x_n))=$$ $$df(x_1,\ldots,x_n).$$}
{\small $$dg(x_2,\ldots,x_n)=$$
$$(-1)^n\sum_{y\in Q}\sum_{i=2}^n(-1)^i((y,\ldots,x_{i-1},x_{i+1},\ldots,x_n)-(y*x_i,\ldots,x_{i-1}*x_i,x_{i+1},\ldots,x_n))=$$
$$(-1)^n\sum_{y\in Q}\sum_{i=3}^n(-1)^i((y,\ldots,x_{i-1},x_{i+1},\ldots,x_n)-(y,\ldots,x_{i-1}*x_i,x_{i+1},\ldots,x_n))=$$
$$g(\sum_{i=3}^n(-1)^{i+1}((x_2,\ldots,x_{i-1},x_{i+1},\ldots,x_n)-(x_2*x_i,\ldots,x_{i-1}*x_i,x_{i+1},\ldots,x_n))=$$
$$gd(x_2,\ldots,x_n).$$}
Furthermore, as shown in
\cite{L-N}(Section 3), 
the composition $gf$ is chain homotopic to $|Q|Id$. Therefore, on homology,
$(gf)_* = |Q|Id$. On the other hand, $(gf)_* = g_*f_*$. Thus,
if $Nf_*(a)=0$, then $0=Ng_*f_*(a)=N|Q|a$, and the theorem follows.
\end{proof}

\section{Appendix 2}

In this section, we work exclusively with quandle homology, and
we show that for $R_k$, $k$ odd, we can compute homology
from a smaller chain complex of ``symmetric" elements.\\
We will work with homology with coefficients in any ring $\mathcal{R}$ (treated as an abelian group).
Let $C^{inv(a)}_n(Q; \mathcal{R})$ denote the subgroup of 
$C_n(Q; \mathcal{R})$ composed of elements $w$ which satisfy $w=w*a$.
$C^{inv(a)}_*(Q; \mathcal{R})$ is a subchain complex of 
$C_*(Q, \mathcal{R})$. Indeed, $*_a(dw)=d(*_aw)=dw$.\\
Let 
$C^{cov(a)}_n(Q; \mathcal{R})=C_n(Q; 
\mathcal{R})/C^{inv(a)}_n(Q; \mathcal{R})$. 
We have a short exact sequence of chain complexes:
$$ 0 \to C^{inv(a)}_*(Q; \mathcal{R}) \stackrel{i}{\to} C_*(Q; \mathcal{R}) 
\to C^{cov(a)}_*(Q; \mathcal{R}) \to 0,$$
and related long exact sequence of homology
$$ \ldots \to H_{n+1}^{cov(a)}(Q; \mathcal{R}) \stackrel{\partial}{\to} 
H_n^{inv(a)}(Q; \mathcal{R})
 \stackrel{i_*}{\to} H_{n}(Q; \mathcal{R})\to$$
$$ \to H_n^{cov(a)}(Q; \mathcal{R}) 
\stackrel{\partial}{\to} 
H_{n-1}^{inv(a)}(Q; \mathcal{R}) \stackrel{i_*}{\to} 
H_{n-1}(Q; \mathcal{R}) \to \ldots .$$

\begin{proposition}\label{h_a_restricted}
The map $h_a$ restricted to $C^{inv(a)}_n(Q; \mathcal{R})$ is a chain map.
\end{proposition}
\begin{proof}
$dh_a(w)=h_ad(w)+(-1)^n(w-w*a)=dh_a(w)$,
as $w\in C^{inv(a)}_n(Q; \mathcal{R})$.
\end{proof}

\begin{lemma}\label{Lemma_split}
If $\mathcal{R}$ is a ring with $2$ invertible, and $Q$ is a kei, then 
the above short exact sequence of chain complexes
splits. In particular, the connecting homomorphism $\partial$, in 
the long exact sequence of homology, is a zero map.
\end{lemma}
\begin{proof} 
Let $\{C^+_n\}$ be a subchain complex of $\{C_n(Q;\mathcal{R})\}$ 
generated by the symmetric elements of the form $u+u*a$, where $u\in Q^n$. 
Similarly, let $\{C^-_n\}$ be a subchain complex of 
$\{C_n(Q;\mathcal{R})\}$ generated by the antisymmetric elements $u-u*a$.
Notice that $*_a(w)=w$ for $w\in C^+_n$, and $*_a(v)=-v$ for $v\in C^-_n$.
It follows that $C^+_n\cap C^-_n=0$. Indeed, if $w\in C^+_n\cap C^-_n$, then 
$*_a(w)=w$ and $*_a(w)= -w$, therefore, $w=0$ ($C^+_n$ and $C^-_n$ are 
$\mathcal{R}$-modules with $2$ invertible in $\mathcal{R})$.\\
Furthermore, $C_n(Q; \mathcal{R})=C^+_n\oplus C^-_n$, because for 
any element $w\in C_n(Q; \mathcal{R})$, 
we have $w+w*a + w- w*a\in C^+_n\oplus C^-_n$,
so $2w\in C^+_n\oplus C^-_n$. 
Since $2$ is invertible in $\mathcal{R}$, $w\in C^+_n\oplus C^-_n$.\\
Thus, the above short exact sequence of chain complexes splits.
From the fact that 
$C^+_n\subset C^{inv(a)}_n(Q; \mathcal{R})$ and  
$C^-_n\cap C^{inv(a)}_n(Q; \mathcal{R})=0$, follows that
$C^+_n=C^{inv(a)}_n(Q; \mathcal{R})$. 
Similarly, $C^-_n=C^{cov(a)}_n(Q; \mathcal{R})$.
\end{proof}

\begin{lemma}
If $Q$ is a kei, and the homomorphism $2Id: H_{n}(Q;\mathcal{R}) \to 
H_{n}(Q;\mathcal{R})$ is an isomorphism (e.g., $H_n$ is a 
finite abelian group without even torsion, or $\Z[\frac{1}{2}]$-module), 
then the homomorphism $i_* : H_n^{inv(a)}(Q; \mathcal{R}) 
\to H_{n}(Q; \mathcal{R})$ is an epimorphism.
\end{lemma}
\begin{proof} 
If the cycles $w_1,\ldots, w_s$ are generators of $H_{n}(Q; \mathcal{R})$,
then also\\ $2w_1,\ldots, 2w_s$ generate this group. By Proposition \ref{Proposition 2.1}(i),
$w+ w*a =2w$ in $H_{n}(Q; \mathcal{R})$, therefore, $a$-symmetric elements 
generate $H_{n}(Q; \mathcal{R})$, and the lemma follows. 
\end{proof} 

As a consequence, we have the result which says that for odd dihedral 
quandles $R_k$, the homology $H_n^{inv(a)}(R_k)$ determines the quandle 
homology $H_n(R_k)$. More precisely: 
\begin{theorem}\label{recover_hom}
For $n>1$, the homomorphism $H_n^{inv(a)}(R_k) \to H_n(R_k)$,
yielded by embedding of chain complexes, is an epimorhism, with the 
kernel being a 2-torsion part of $H_n^{inv(a)}(R_k)$.
\end{theorem}
\begin{proof}
We use Lemmas 8.2 and 8.3 for homology of $R_k$, and the ring 
$\mathcal{R}=\Z[\frac{1}{2}]$. We get the isomorphism: 
$$ H_n^{inv(a)}(R_k;\Z[\frac{1}{2}]) \to 
H_n(R_k;\Z[\frac{1}{2}]).$$ 
From the universal coefficient theorem for homology, and the fact that 
$\Z[\frac{1}{2}]$ is a torsion free group, we get the isomorphism
$$ H_n^{inv(a)}(R_k)\otimes \Z[\frac{1}{2}] \to
H_n(R_k) \otimes \Z[\frac{1}{2}].$$ The theorem follows, 
because $H_n(R_k)$ is a group without a $2$-torsion \cite{L-N}.
\end{proof}

\begin{center}
\begin{tabular}{l  @{\hspace{2.5 cm}} l}
                            &                         \\
Maciej Niebrzydowski        &     J\'ozef~H. Przytycki\\
e-mail: mniebrz@gmail.com      &     e-mail: przytyck@gwu.edu\\
                            &
\end{tabular}
\end{center}
J\'ozef Przytycki is on sabbatical leave at
University of Maryland, College Park.

\vspace{0.4 cm}

\noindent \textsc{Dept. of Mathematics, Old Main Bldg., 1922 F St. NW\\
The George Washington University, Washington DC, 20052}
\end{document}